\documentstyle[amsfonts,12pt]{article}
\def\qed{\hfill \rule{.5em}{1em} \medskip}
\newtheorem{example}{Example}
\newtheorem{definition}{Definition}
\newtheorem{theorem}{Theorem}

\newtheorem{lemma}{Lemma}
\def\ole{ \overline}
\begin{document}

\centerline{\Large{\bf K\"ahler Yamabe minimizers on minimal ruled surfaces}}

\vspace{10mm}

\centerline{\bf Christina W. T{\o}nnesen-Friedman}
\centerline{Department of Mathematical Sciences}
\centerline{University of Aarhus, Denmark}

\vspace{8mm}

\begin{abstract} It is shown that if a minimal ruled surface ${\mathbb P}(E) 
\rightarrow \Sigma$ admits a K\"ahler Yamabe minimizer then this
metric is generalized K\"ahler-Einstein and the holomorphic vector bundle
$E$ is quasi-stable.
\end{abstract}

\section{Introduction}

The {\it minimal ruled surfaces} form a special class of compact K\"ahlerian
surfaces and are by definition the
total spaces of ${\mathbb CP}_1$ bundles over compact Riemann surfaces $\Sigma$.
Any ruled surface can be written \cite{bbv} as
\[
{\mathbb P}(E) \rightarrow \Sigma
\]
i.e., as the projectivization of a holomorphic
rank two vector bundle $E$ over $\Sigma$, where $E$ is unique up to 
tensoring with
a holomorphic line bundle. Moreover any ruled
surface is birationally equivalent to $\Sigma \times {\mathbb CP}_1$. In 
particular, any ruled surface is algebraic. In fact, the minimal
models of any complex surface which is birationally equivalent to $\Sigma
\times {\mathbb CP}_1$, are exactly the ruled surfaces \cite{beau,yang}.

Suppose that $E \rightarrow \Sigma$ is {\it quasi-stable}, that is, $E$ is 
semi-stable
(in the sense of Mumford) and decomposes into a direct sum 
\[
E = E_1 \oplus...\oplus E_k\]
of stable sub-bundles (here $k=1$ or $2$) such that 
\[
\frac{deg(E)}{rank(E)} = \frac{deg(E_i)}{rank(E_i)}
\]
for $i=1,...,k$. 
Narasimhan and Seshadri \cite{28} have proved that quasi-stability is 
equivalent to the existence of a flat projective unitary connection on $E$.
In other words, if $E$ is quasi-stable, then ${\mathbb P}(E) \rightarrow \Sigma$ 
is a flat ${\mathbb CP}_1$ bundle,
i.e., is defined by some representation 
\[
\rho : \pi_1(\Sigma) \rightarrow {\mathbb P}SU(2) = SO(3).
\] 
So, when $E$ is quasi-stable, local products of constant scalar curvature
K\"ahler metrics on $\Sigma$ and ${\mathbb CP}_1$ will exhaust the entire K\"ahler
cone on the ruled surface with K\"ahler classes of constant scalar curvature
K\"ahler metrics.

Burns and de Bartolomeis proved that quasi-stability is a necessary condition 
for 
the existence of scalar-flat
K\"ahler metrics. More recently 
LeBrun 
proved a similar statement for negative constant scalar curvature.
We summarize these results in the theorem below.

\begin{theorem}
{\bf (Burns, de Bartolomeis \cite{29} and LeBrun \cite{leb2})}
Let ${\mathbb P}(E) \rightarrow \Sigma$ be a minimal ruled surface with a 
K\"ahler class $[\omega]$ such that $c_1 \cdot [\omega] \leq 0$. Then 
$[\omega]$ contains
a K\"ahler metric of constant scalar curvature if and only if $E \rightarrow
\Sigma$ is a quasi-stable vector bundle.
\label{bbl}
\end{theorem} 

A key step in both proofs is the observation that the constant scalar
curvature K\"ahler metric must be K\"ahler with respect to two non-equivalent
complex structures on the ruled surface.

Whether the statement holds in the case $c_1 \cdot [\omega] >0$ is still 
unknown.
In this paper we assume that the K\"ahler metric is also a Yamabe
minimizer in its conformal class and show that then quasi-stability holds.

\section{Perturbed Seiberg-Witten Invariants}

Let $M$ be a compact, oriented four manifold such that $H^2(M,{\mathbb R})$ has
dimension two and $b_+ = b_- = 1$. (In general, one could let $b_-$ have 
arbitrary value.) Let $g$ be a Riemannian metric on $M$ and $\star$ be the
Hodge Star operator defined with respect to $g$ and the orientation. 
Then the
one dimensional subspace of $H^2(M,{\mathbb R})$
\[
H^+(g) := \{ [\nu] \in H^2(M,{\mathbb R}) \mid \star \nu = \nu \}
\]
is called a {\it metric polarization} \cite{leb2}. Observe that
\[
H^-(g) := \{ [\nu] \in H^2(M,{\mathbb R}) \mid \star \nu =- \nu \} 
\]
is the metric polarization with respect to the opposite orientation. If $g$
is K\"ahler then the K\"ahler class spans $H^+(g)$.

The open cone
\[
\{ [\nu] \in H^2(M,{\mathbb R}) \mid [\nu] \cdot [\nu] > 0 \}
\]
consists of two connected components, called {\it nappes} \cite{leb3}. Given 
a nappe
${\cal C}^+$ and a Riemannian metric $g$, let $\omega$ be a $g$-harmonic,
self-dual two form such that $[\omega] \in {\cal C}^+$. This form always
exists and is unique 
up to multiplication with a positive constant. In fact $[\omega] \in H^+(g)
\cap {\cal C}^+$.
If $M$ has a K\"ahler metric then the canonical choice of nappe is the
one containing the K\"ahler class. This way the corresponding 
$g$-harmonic, self-dual two-form 
to any metric on $M$ is simply a generalization of the K\"ahler form. 

Now assume that $M$ has a $Spin_c$ structure $c$ of almost-complex type. 
Relative to any metric $g$, the {\it perturbed Seiberg-Witten
invariant}
$p_c(M,{\cal C}^+)$ is defined to be the number of solutions, modulo gauge 
and counted with
orientations, of the perturbed Seiberg-Witten equations \cite{krm,wit}
\begin{equation}
D_A \Phi = 0
\label{1}
\end{equation}
\begin{equation}
iF^+_A + \sigma(\Phi) = \epsilon
\label{2}
\end{equation}
where $\epsilon$ is a generic (so that $(g,\epsilon)$ is excellent) 
self-dual two
form with $\int_M \epsilon \wedge \omega > 2\pi c_1(c) \cdot [\omega]$.
Note that all $\epsilon$ satisfying this inequality make $(g,\epsilon)$ a
good pair and it is easy to see that they are all in the 
same chamber. The above invariant is therefore well-defined and metric 
independent. 
We refer to \cite{leb1} for definitions of the words ``excellent'', ``good'' 
and ``chamber''.

Note that if $p_c(M,{\cal C}^+) \neq 0$ the the equations (\ref{1}) and 
(\ref{2}) have a solution $\Phi \neq 0$ for any $\epsilon = t\omega$
where $t \gg 0$. This is easily seen by the fact that $(g,t\omega)$ is a good
pair (in the chamber determined by ${\cal C}^+$). If it had no solutions it 
would automatically
be excellent and therefore contradict the non-vanishing of $p_c(M,{\cal C})$.
Therefore, $(g,t\omega)$ has a solution (not necessarily transverse) and by
$(g, t\omega)$ being good this solution is irreducible ($\Phi \neq 0$).

\begin{example}\cite{leb1}
If $(M,J,g)$ is a K\"ahler surface, $c$ the $Spin_c$ structure induced by
$J$ and ${\cal C}^+$ the canonical choice of nappe then
$p_c(M,{\cal C}^+) \neq 0$.
\label{e2}
\end{example}

The perturbed Seiberg-Witten invariant is also defined for $Spin_c$ structures
on $M$ who do not arise from an almost-complex structure \cite{LiLiu}.

If $(M,g,J)$ is an almost-K\"ahler manifold then the almost-K\"ahler form 
$\omega$
is a harmonic self-dual form. Hence, even though $J$ may not be integrable,
$[\omega]$ still determines a canonical
choice of nappe. Since $\vert \omega \vert = \sqrt{2}$
the following result is a direct application of (\cite{leb3}, Theorem 1 and
2).

\begin{theorem}
{\bf (LeBrun)}
Let $(M,g,J)$ be an almost-K\"ahler surface with the canonical choice of
nappe ${\cal C}^+$. If $c$ is a $Spin_c$ structure such that $p_c({\cal C}^+,M)
\neq 0$ then 
\[
\int_M s d\mu \leq 4\pi c_1(c) \cdot [\omega]
\]
where $s$ is the scalar curvature, $d\mu$ is the metric volume form and 
$c_1(c) = c_1(\det V^+)$. Moreover, 
equality is achieved if and only if $(M,g,J)$ is K\"ahler and $J$ is 
compatible with $c$.
\label{hurra}
\end{theorem}
For the proof we refer to LeBrun's paper \cite{leb3}. However in (\cite{leb3},
Theorem 2)
the compability statement was made without offering a proof. For the sake
of completeness we now prove this. 
When equality is achieved
we have that $(M,g,J)$ is K\"ahler. Therefore
\[
4\pi c_1(c) \cdot [\omega] = \int_M s d\mu = 4\pi c_1(K^{-1}) \cdot [\omega]
\]
where $K^{-1}$ is the anti-canonical line bundle of $(M,J)$. The
compability of $J$ with $c$ then follows from the following lemma.
\begin{lemma}
Let $M$ be a compact smooth manifold with $b_+=1$. Assume that $M$ has a 
K\"ahler metric $g$ with K\"ahler form $\omega$ and complex structure $J$.
Let $K$ denote the canonical line bundle of $(M,J)$.
Let ${\cal C}^+$ be the canonical choice of nappe.
Suppose $c$ is any $Spin_c$ structure on $M$ with corresponding complex 
line bundle
$L=\det V^+$ such that $p_c(M,{\cal C}^+) \neq 0$.
Then $E=(K \otimes L)^{\frac{1}{2}}$ is either trivial or a holomorphic 
line bundle corresponding
to an effective divisor. In particular, $c_1(L) \cdot [\omega] \geq c_1(K^{-1})
\cdot [\omega]$ with equality if and only if $E$ is trivial and 
$c$ is the $Spin_c$ structure induced by $J$.
\label{tmf}
\end{lemma}

\noindent{\bf Proof:}
The trick is to choose the perturbation to be $\epsilon = t\omega$, $t \gg 0$. 
Now we follow Witten's calculations for the unperturbed Seiberg-Witten
equations on a K\"ahler manifold \cite{wit} (see also the proof of
(\cite{mf}, Proposition 2.1.)).
Since
$\epsilon$ is of type $(1,1)$ with respect to the K\"ahler structure, we get
by precisely the same argument as in \cite{wit} that for a solution
$(A,\Phi)$ to both equation (\ref{1}) and
\begin{equation}
iF_A^+ + \sigma(\Phi) = t\omega
\label{3}
\end{equation}
the curvature $F_A$ is of type $(1,1)$ and $E$ has a holomorphic structure
(induced by $D_A$). If we write $\Phi = (\alpha,\beta)$ where $\alpha$ is
a section of $E$ and $\beta$ is a section of $\bigwedge^{0,2}(E)$ then
$\alpha$ and $\ole{\beta}$ are holomorphic and one of them must vanish.
Now (\ref{3}) rewrites to
\[
iF_A^+ = \frac{(-\vert \alpha \vert^2 + \vert \beta \vert^2 + 4t)}{4} \omega
\]
implying that
\[ 
2\pi c_1(L) \cdot [\omega] = \frac{(-\vert \alpha \vert^2 + \vert \beta 
\vert^2 + 4t)}{4} [\omega]^2.
\]
For $t$ sufficiently large we must have that $\alpha$ is a non-vanishing 
holomorphic section of $E$. Thus, unless it is trivial,
the line bundle $E$ corresponds to an effective divisor.
The inequality now follows from the fact that
the ``area'' of any effective divisor on the K\"ahler manifold is non-zero.

If $E$ is trivial then $L=K^{-1}$ and, since a $Spin_c$ structure on an
almost-complex manifold is determined by the determinant line bundle 
$L=\det V^+$, we are done.

\qed

The author would like to point out that Lemma \ref{tmf} a special case of 
Theorem 1.3. in \cite{taubes} where Taubes
proved a similar statement in the symplectic setting.

\section{The Yamabe Constant}

\begin{definition}
Let $g$ be a Riemannian metric on a four manifold $M$. The {\it Yamabe
constant} of the corresponding conformal class $[g]$ is defined to be
\[
Y_{[g]} = \inf_{g \in [g]} \frac{\int_M s_g d\mu_g}{(\int_M d\mu_g)^
{\frac{1}{2}}}
\]
\end{definition}
Note that the above infimum is in fact achieved by a metric in $[g]$. This
was proved by Yamabe, Trudinger, Aubin
and Schoen \cite{aub,lp,sch}. A metric which minimizes $\frac{\int_M s_g d\mu_g}{(\int_M d\mu_g)^{\frac{1}{2}}}$ on $g$ is called a {\it Yamabe
minimizer}. Any Yamabe minimizer must have constant scalar curvature.
If $Y_{[g]} \leq 0$ then $g$ is the unique (up to scalar multiplication)
Yamabe minimizer of $[g]$ if and only if $g$ has constant scalar curvature.
Unfortunately, for $Y_{[g]} > 0$ constant scalar curvature does not necessarily
imply that a metric is a minimizer and uniqueness of the minimizers does
not always hold in this situation either. Observe that $Y_{[g]} >0$ if and
only if there exists a metric in $[g]$ with strictly positive scalar curvature.

Applying Theorem \ref{hurra}, LeBrun found an estimate for $Y_{[g]}$.

\begin{theorem}
{\bf (LeBrun\cite{leb3})}
Let $(M,[g])$ be an oriented conformal Riemannian four-manifold, and let
$\omega$ be a closed $2$-form which is self-dual with respect to $[g]$ and not
identically zero. Suppose that $b^+(M) =1$ and that the perturbed 
Seiberg-Witten 
invariant $p_c(M,{\cal C}^+)$ is non-zero for some $Spin_c$ structure
$c$, where ${\cal C}^+ \subset H^2(M,{\mathbb R})$ is the nappe containing 
$[\omega]$.
Then the Yamabe constant of $[g]$ satisfies
\[
Y_{[g]} \leq \frac{4\pi c_1(c) \cdot [\omega]}{\sqrt{[\omega]^2/2}}.
\]
Moreover, equality is achieved if and only if there is a Yamabe minimizer $g 
\in [g]$
which is K\"ahler, with K\"ahler form $\omega$ and complex structure compatible
with $c$.

\label{yamtheo}
\end{theorem}

\begin{definition}{\cite{mat}}
A K\"ahler metric is said to be generalized K\"ahler-Einstein if the Ricci
form is parallel with respect to the Levi-Civita connection.
\end{definition}

We can now prove the following theorem.
\begin{theorem}
Let $M = {\mathbb P}(E) \rightarrow \Sigma$ be a minimal ruled surface over a 
compact Riemann surface $\Sigma$. If $M$ has a K\"ahler metric $g$ with
constant positive scalar curvature such that $g$ is a
Yamabe minimizer in $[g]$, then $g$ is generalized K\"ahler-Einstein
and therefore locally a product. Consequently $E$ is a quasi-stable 
holomorphic vector bundle.
\label{p1}
\end{theorem}
If $g$ has constant non-positive scalar curvature then the above is true by
Theorem \ref{bbl}.

\noindent{\bf Proof:}
First assume that $\Sigma$ is ${\mathbb CP}_1$. The only
Hirzebruch surface with constant scalar curvature K\"ahler metric is
the product ${\mathbb CP}_1 \times {\mathbb CP}_1$ \cite{cal}. On this 
surface any constant scalar curvature K\"ahler metric must be invariant
under the $SO(3)$ action on each ${\mathbb CP}_1$. This forces the
metric to be a product of (multipla of) the 
Fubini-Study metric. We have used the fact that any extremal K\"ahler
metric is invariant under the action of the 
maximal compact subgroup
of (the identity component of) the group of holomorphic 
transformations\cite{cal2}.

Now assume that the genus ${\bf g}$ of
$\Sigma$ is at least one. Let $g$ be a K\"ahler Yamabe minimizer with
positive scalar curvature. The Yamabe constant of $[g]$ is then given by
\begin{equation}
Y_{[g]} = \frac{4\pi c_1 \cdot [\omega]}{\sqrt{[\omega]^2/2}}
\label{equality}
\end{equation}
where $\omega$ is the K\"ahler form of $g$ and $c_1 = c_1(K^{-1})$.
Note that for a minimal ruled surface $b_+ = b_- = 1$. 
Let $c$ be the $Spin_c$ structure induced by the complex structure $J$ on $M$
and let ${\cal C}^+$ be the canonical choice of nappe. According to 
Example \ref{e2}, $p_c(M,{\cal C}^+) \neq 0$ and 
we see that equation (\ref{equality}) is a special case of 
Theorem \ref{yamtheo}.

Now consider the fiber-wise anti-podal map 
$\psi : M \rightarrow M$ \cite{leb2}.
This is an orientation reversing diffeomorphism
and we can define a $Spin_c$ structure $\ole{c}$
on $\ole{M}$ by setting $\ole{c} := \psi^\ast c$.
Observe that $\psi^\ast$ sends ${\cal C}^+$ to a nappe $\psi^\ast {\cal C}^+$
for $\ole{M}$ and $\psi^\ast(H^+(g)) = H^-(g)$ \cite{leb2}.
Since $(\ole{M},\ole{c},\psi^\ast{\cal C}^+)$ and $(M,c,{\cal C}^+)$ are
isomorphic as oriented four-manifolds with nappes and $Spin_c$ structures
we have that
\[
p_{\ole{c}}(\ole{M},\psi^\ast{\cal C}^+) = p_c(M,{\cal C}^+) \neq 0.
\]
Theorem \ref{yamtheo} applied to $(\ole{M},\ole{c}, \psi^\ast{\cal C}^+)$ 
now implies that
\begin{equation}
Y_{[g]} \leq \frac{4\pi \psi^\ast c_1 \cdot \psi^\ast [\omega]}
{\sqrt{(\psi^\ast[\omega])^2/2}}
\label{ineq}
\end{equation}
on $\ole{M}$. But the Yamabe constant is independent of orientation and the
right hand side of (\ref{ineq}) is just the right hand side of 
(\ref{equality}). So we must have equality in (\ref{ineq}) and thus there
exist a Yamabe minimizer $\tilde{g} \in [g]$ such that $\tilde{g}$ is 
K\"ahler with
respect to some complex structure $\tilde{J}$ in $\ole{c}$ where the K\"ahler
form $\tilde{\omega}$ is equal to the harmonic part of $\psi^\ast \omega$. 

Now we want to show that $\tilde{g}=g$. We can
assume that $\int d\mu = \int d\tilde{\mu} = 1$. If we write $\tilde{g} =
u^2 g$ for some positive smooth function $u$ we have that 
\begin{equation}
\int u^4 d\mu =1
\label{intu4}
\end{equation}
and
\[
\vert \tilde{\omega} \vert^2 = u^4\tilde{\vert}\tilde{\omega}\tilde{\vert}^2
=2u^4.
\]
Since $\tilde{s}=Y_{[g]}=s$ we have that
\begin{equation}
\Delta u = \frac{s(u^3-u)}{6}.
\label{Deltau}
\end{equation}

Since the Euler characteristic of $M$ is given by $\chi = 4(1-{\bf g})$ and the
signature $\sigma = b_+ - b_-$ vanishes, the (strict) Hitchin-Thorpe 
inequality \cite{hitchin,thorpe,besse}, $2\chi > 3\vert \sigma \vert$,
is not satisfied when ${\bf g}>0$. Therefore no Riemannian metric on 
$M$ can be Einstein.
In particular, the primitive part $\rho_0$ (resp. $\tilde{\rho_0}$) of the 
Ricci form $\rho$ (resp. $\tilde{\rho}$) of $g$ (resp. $\tilde{g}$) does 
not vanish. Moreover,
$d\rho_0 = d\tilde{\rho_0} =0$, which follows from the fact that the scalar
curvatures are constant. 
Since $b_+=b_-=1$ we must therefore have that $\tilde
{\omega} = k\rho_0$ and $\omega = \tilde{k}\tilde{\rho_0}$
where $k$ and $\tilde{k}$ are non-zero constants.
Now $\psi^\ast[\frac{\rho}{2\pi}]=\psi^\ast c_1= c_1(\ole{c})=
[\frac{\tilde{\rho}}{2\pi}]$. In particular 
\[
\psi^\ast[\rho_0]=[\tilde{\rho_0}]
\]
thus
\[
k^{-1}\psi^\ast[\tilde{\omega}] = \tilde{k}^{-1} [\omega]
\]
hence
\[
k^{-1}\psi^\ast\psi^\ast[\omega] = k^{-1}[\omega]=\tilde{k}^{-1} [\omega]
\]
and consequently $k=\tilde{k}$.

We can calculate $k$ up to a sign as follows:
\[
\begin{array}{rcl}
c_1^2 & = & \frac{1}{(2\pi)^2}\int \rho \wedge \rho \\
\\
      & = & \frac{1}{(2\pi)^2}(\int (\frac{s}{4}\omega)\wedge 
(\frac{s}{4}\omega) + \int \rho_0 \wedge \rho_0) \\
\\
      & = & \frac{1}{(2\pi)^2}(\frac{s^2}{8} - \int \vert \rho_0 \vert^2 d\mu)
\\
\\
      & = & \frac{1}{(2\pi)^2}(\frac{s^2}{8} - k^{-2}\int \vert \tilde{\omega}
\vert^2 d\mu) \\
\\
      & = & \frac{1}{(2\pi)^2}(\frac{s^2}{8} -2k^{-2}\int u^4 d\mu)\\
\\
      & = & \frac{1}{(2\pi)^2}(\frac{s^2}{8} -2k^{-2})
\end{array}
\] 
and therefore
\begin{equation}
k^{-2} = \frac{s^2}{16} - 2\pi^2 c_1^2 = \frac{s^2}{16} - 2\pi^2(2\chi +
3\sigma)= 
\frac{s^2}{16} + 16\pi^2({\bf g}-1).
\label{k}
\end{equation}
The traceless part of the Ricci tensor of $\tilde{g}$ can now be found
as follows:
\[
\begin{array}{rrr}
\tilde{r}_0(X,Y) & = & \tilde{\rho}_0(X,\tilde{J}Y) \\
\\ 
                 & = & k^{-1}\omega(X,\tilde{J}Y) \\
\\
                 & = & k^{-1}g(JX,\tilde{J}Y)\\
\\
                 & = & k^{-1}u^{-2}\tilde{g}(JX,\tilde{J}Y)\\
\\
                 & = & -k^{-1}u^{-2}\tilde{\omega}(JX,Y) \\
\\
                 & = & -u^{-2}\rho_0(JX,Y) \\
\\
                 & = & u^{-2}r_0(X,Y).
\end{array}
\]
On the other hand, since $\tilde{g}=u^2 g$ we have from (\cite{besse}, 
(1.161b)) that
\[
\tilde{r}_0 = r_0 + 2u(\nabla d(u^{-1}) + \frac{\Delta(u^{-1})}{4} g)
\]
and hence from the above calculation
\[
u^{-2}r_0 = r_0 + 2u(\nabla d(u^{-1}) + \frac{\Delta(u^{-1})}{4} g).
\]
Using equation (\ref{Deltau}) we find that
\[
\begin{array}{rrr}
\Delta(u^{-1}) & = & -2u^{-3}\vert du \vert^2 - u^{-2}\Delta u\\
               & = & -2u^{-3}\vert du \vert^2 -\frac{s(u-u^{-1})}{6}
\end{array}
\]
and therefore
\[
\nabla d(u^{-1}) = \frac{(u^{-3} - u^{-1})}{2} r_0 + (\frac{u^{-3}\vert du 
\vert^2}{2} +  \frac{s(u-u^{-1})}{24})g.
\]
In particular, at a maximum of $u^{-1}$ the Hessian of $u^{-1}$ is
given by
\begin{equation}
\nabla d(u^{-1}) = \frac{(u-u^{-1})}{2}(\frac{s}{12}g - u^{-2}r_0).
\label{hessiant}
\end{equation}

Let $p \in M$ be any point on our manifold.
Since $r_0$ is a traceless symmetric tensor of type $(1,1)$
we can find an orthonormal base $\{ e_1,Je_1,e_2,Je_2 \}$ of $T_pM$ such that
$r_0$ can be represented by the matrix
\[
\left(
\begin{array}{cccc}
\lambda           & 0       & 0        & 0 \\
0                 & \lambda & 0        & 0 \\
0                 & 0       & -\lambda & 0 \\
0                 & 0       & 0        & -\lambda
\end{array}
\right)
\]
where $\lambda$ is the positive constant given by
\[
\lambda = \frac{\vert r_0 \vert}{2} = \frac{ \vert \rho_0 \vert }{\sqrt{2}} =
\frac{\vert \tilde{\omega} \vert}{\vert k \vert \sqrt{2}} = \frac{u^2}
{\vert k \vert}.
\]
The tensor $(\frac{s}{12}g-u^{-2}r_o)$ can now be represented by the matrix
\[
\left(
\begin{array}{cccc}
\frac{s}{12} - \frac{1}{\vert k \vert} & 0 & 0 & 0\\
0 & \frac{s}{12} - \frac{1}{\vert k \vert} & 0 & 0\\
0 & 0 & \frac{s}{12} + \frac{1}{\vert k \vert} & 0 \\
0 & 0 & 0 & \frac{s}{12} + \frac{1}{\vert k \vert} 
\end{array}
\right)
\]
and since
\[
\begin{array}{ccl}
\frac{s}{12} - \frac{1}{\vert k \vert} & =    & \frac{s}{12} - \sqrt{\frac{s^2}
{16} + 16\pi^2({\bf g}-1)}\\
                                       & \leq & \frac{s}{12} - \sqrt{\frac{s^2}
{16}}\\
                                       & = & -\frac{s}{6}\\
                                       & < & 0
\end{array}
\]
the tensor is never semi-definite.
But at the maximum of $u^{-1}$ the Hessian must be negative semi-definite and
hence from equation (\ref{hessiant}) we have that $u^{-1} = u = 1$ at
the maximum of $u^{-1}$ and by equation (\ref{intu4}) we conclude that $u=1$
and hence $\tilde{g}=g$ everywhere.

Now the Ricci form satisfies 
\[
\nabla \rho = \frac{s}{4}\nabla\omega + \nabla \rho_0 =0.
\] 
Thus $g$ is generalized K\"ahler-Einstein. Since $\tilde{g}=g$
(or since $g$ is generalized K\"ahler-Einstein with
non-vanishing $\rho_0$) we have that $g$ is K\"ahler
with respect to two complex structures $J$ and $\tilde{J}$ inducing opposite 
orientations. Therefore
the holomony \cite{besse} is a subgroup of $U(1) \times U(1)$ and the universal
cover $(\hat{M}, \hat{g})$ of $(M,g)$ must be a Riemannian product
$(\hat{M}, \hat{g}) = (M_1,g_1) \times (M_2, g_2)$ of a pair of 
complete simply connected surfaces. Clearly the scalar curvature of each
$(M_i,g_i)$ must be constant and
since $s>0$ (but also for topological reasons \cite{leb2}) we must have that 
at least one of the surfaces is a two sphere. Thus 
$(\hat{M}, \hat{g}) = S^2 \times (M_2, g_2)$.
Since the genus of $\Sigma$ is at least one
$(M_2,g_2)$ must be either ${\mathbb C}$ or ${\mathbb CH}_1$ with their
canonical metric.
The rest of the proof follows along the same line of reasoning as in
the proof of
(\cite{leb2}, Theorem 4). In order to make this paper reasonably 
self-contained we repeat the arguments here.  
The holomony of $(\hat{M}, \hat{g})$ is $U(1) \times U(1)$ so 
the lift of $J$ on $M$ must coincide with the
product complex structure, once the factors are correctly oriented.
Since the genus of
$\Sigma$ is at least one we have that $\pi_1(M) = \pi_1(\Sigma)$ is
non-trivial and acts on $S^2 \times (M_2, g_2)$ by holomorphic isometries
sending any compact holomorphic curve $S^2 \times \{ pt \}$ to another
curve of this form. The induced action on $M_2$ is moreover free and proper,
since $S^2$ is compact and every rotation of $S^2$ has a fixed point. Thus
$M = (S^2 \times M_2)/\pi_1(\Sigma)$ is biholomorphic to $\tilde{\Sigma}
\times_\rho {\mathbb CP}_1$ for some compact Riemann surface $\tilde{\Sigma}$ 
and some representation $\rho : \pi_1(\tilde{\Sigma}) \rightarrow 
{\mathbb P}SU(2)=SO(3)$. By uniqueness of ruling this biholomorphism must
be a bundle biholomorphism inducing a biholomorphism between
$\tilde{\Sigma}$ and $\Sigma$.
Thus $M = {\mathbb P}(E) \rightarrow \Sigma$ is a flat ${\mathbb CP}_1$ 
bundle and $E$ is therefore quasi-stable.
\qed

\noindent {\bf Acknowledgment}
The author would like to thank C. LeBrun, J.E. Andersen, K. Akutagawa
and H. Pedersen for very helpful conversations.

\end{document}